%% file: thesis.tex
\documentclass{article}
\usepackage{amsmath,amsfonts,epsfig}
\input{macros}

\newcommand{\dev}{D}

\renewcommand{\hbar}{\bar{{\mathbb H}}^3}

\begin{document}

\title{Rigidity of geometrically finite hyperbolic cone-manifolds}

\author{K. Bromberg\footnote{This work was partially supported by grants from the NSF}}

\date{October 29, 2002}

\maketitle

\begin{abstract}
\noindent
In a recent paper Hodgson and Kerckhoff
\cite{Hodgson:Kerckhoff:cone} prove a local
rigidity theorem for finite volume, 3 dimensional hyperbolic
cone-manifolds. In this paper we extend this result to geometrically
finite cone-manifolds. Our methods also give a new proof of a local
version of the classical rigidity theorem for geometrically finite
hyperbolic 3-manifolds.
\end{abstract}

\section{Introduction}

A hyperbolic cone-manifold is a singular hyperbolic structure where
the singularity is a simple closed curve with cross section a
hyperbolic cone. We say a hyperbolic structure on a manifold, $M_0$, is
{\em locally rigid} if for any
smooth family of hyperbolic metrics $M_t$, $M_0$ is isometric to
$M_t$ for small $t$.
In a recent paper \cite{Hodgson:Kerckhoff:cone} Hodgson and Kerckhoff
prove a local rigidity result for finite volume 3-dimensional
hyperbolic cone-manifolds. In this paper we extend this result to
geometrically finite cone-manifolds without rank one cusps. The
methods employed were first developed by Calabi
\cite{Calabi:rigidity} and Weil \cite{Weil:compact} in their proof
that closed hyperbolic manifolds of dimension $\geq 3$ are locally
rigid. Garland \cite{Garland:rigidity} extended their result to finite
volume hyperbolic manifolds.

The result here is another example of the rich deformation theory of
hyperbolic manifolds that is special to dimension three.
The contrast between 3-dimensional hyperbolic manifolds and those of
dimension $\geq 4$ can be seen in Garland and Raghunathan's
\cite{Garland:Raghunathan:rigidity} proof that
finite volume hyperbolic manifolds of dimension $\geq 4$ cannot be
deformed even through incomplete hyperbolic structures while
in dimension 3 Thurston \cite{Thurston:book:GTTM} showed that if
$M$ has non-empty boundary there is at least a 1-dimensional 
space of deformations through incomplete structures. The basic
philosophy is that a hyperbolic structure is determined by its
boundary geometry. We will exploit Thurston's result in an essential
way here.

A geometrically finite hyperbolic structure on the interior of a 3-manifold $M$ extends to a conformal structure on the boundary of $M$. Our work here provides a new proof of the 
following well known result which is  the work of
many people including Ahlfors, Bers, Kra, Marden, Maskit, Mostow and
Prasad. An expository account can be found in \cite{Bers:def}.

\begin{theorem}
\label{localnocone}
$M$ is locally rigid rel the conformal boundary.
\end{theorem}

If $\del M$ is incompressible then the classical proof shows that this
is a global result. In the general case there is also a global result
although it takes more work to state.

Although the result is stronger than the local
theorem we prove the methods of proof  of the global theorem do not generalize to 
cone-manifolds. In the original proof one uses the completeness of the
hyperbolic structure to convert the problem to that of studying the
action of a discrete group of M\"obius transformations on the Riemann
sphere, $\chat$. For a cone-manifold the group will not be
discrete and the relationship between the action of the group on
$\chat$ and the hyperbolic structure is unclear. In our proof of
Theorem \ref{localnocone}, little use is made of completeness nor is the
action of the group on $\chat$ studied. In particular, two results at
the heart of the classical proof, the measurable Riemann mapping
theorem and the zero area theorem for limit sets, are not used.

These methods allows us to extend Theorem \ref{localnocone} to
cone-manifolds.

\bigskip

\noindent
{\bf Theorem \ref{main}} 
{\em If $M$ is a geometrically finite cone-manifold without rank one
cusps and all cone angles are $\leq 2 \pi$ then $M$ is locally rigid
rel cone angles and the conformal boundary.}

\bigskip

We remark that Theorem \ref{main} should still hold for structures
with rank one cusps. In particular the classical proof of Theorem
\ref{localnocone} does allow such cusps.
%
%

McMullen \cite{McMullen:graft} has shown that
local rigidity of 
geometrically finite cone-manifolds for cone angles {\em greater} than
$2 \pi$ implies the grafting conjecture for simple closed
curves. This was one the original motivations of this work. Scannell
and Wolf \cite{Scannell:Wolf:graft} have recently 
proved this conjecture for all laminations using harmonic maps. 
Local rigidity for cone angles greater than
$2 \pi$ is still an open question.

The results in this paper also have applications to classical conjectures about Kleinian groups. See \cite{Bromberg:projective}, \cite{Bromberg:schwarz}, \cite{Brock:Bromberg:density} and \cite{Brock:Bromberg:Evans:Souto}.

We now outline the contents of the paper.

The main object of study is the bundle, $E \rightarrow M$, of germs of
Killing fields over the hyperbolic manifold $M$. Weil
\cite{Weil:cohomology} showed that the deRham cohomology group
$H^1(M;E)$ is canonically isomorphic to the Zariski tangent space of
$R(M)$, the space of representations of $\pi_1(M)$ in $PSL_2\cx$
modulo conjugacy. $E$ has a flat connection which gives a covariant
derivative, $d$, and a natural Riemannian metric which allows us to
define a co-derivative, $\delta$, and a Laplacian, $\Delta$.

If $M$ is closed manifold then the Hodge theorem implies that every
cohomology class in $H^1(M;E)$ has a harmonic representative. One then
shows via a Weitzenbock formula that any harmonic representative
is trivial. If the manifold is not closed then we no longer have the
Hodge theorem and a boundary term appears in the Weitzenbock
formula. In this case we prove a Hodge theorem for every cohomology
class that has a 
representative, $\omega$, that is ``almost'' harmonic, in the sense
that $\delta \omega$ is in $L^2$. Our first step is to construct a
model deformation that is ``almost'' harmonic to which we apply the
Hodge	
theorem. We then find an exhaustion of the manifold by 
compact submanifolds and apply the Weitzenbock formula to the harmonic
representative restricted to these compact submanifolds. The last step
is to show the boundary term in the Weitzenbock term approaches zero
as we exhaust the manifold. This will only happen if the cohomology
class represents a deformation that fixes the cone angles and the
conformal boundary.

In  \S\ref{background} we summarize the necessary background
contained in \S1 and \S2 of \cite{Hodgson:Kerckhoff:cone}. We
emphasize those calculations which will be used later in this paper.

\S\ref{ends} is the heart of the paper. In it we construct the
model deformation on a geometrically finite end.

In \S\ref{hodge theory} we prove the Hodge theorem for the model
deformation. We use the Hodge theorem to prove a vanishing theorem for
those cohomology classes that fix the cone angle and the conformal
boundary.

In the conclusion of the paper, \S\ref{representations}, we make
the identification of $H^1(M;E)$ with the Zariski tangent space of
$R(M)$. Using our calculation of $H^1(M;E)$ we show that $R(M)$ is
locally parameterized by the Teichm\"uller space of the conformal boundary
and the complex length of the cone
singularity. This then implies the main result, Theorem \ref{main}.

{\bf Acknowledgments.} The results in this paper are from the
author's thesis. He would like to thank his advisor, Steve Kerckhoff,
for his help throughout this project. He would also like to thank Curt
McMullen for many discussions about this work.
\section{Background deformation theory}
\label{background}


Let $M$ be a manifold and $\rho:\pi_1(M) \longrightarrow PSL_2\cx$ a
representation of its fundamental group. Let $\tilde{E}(M) = \tilde{M} \times sl_2\cx$ and 
let $E(M)$ be the quotient $\tilde{E}(M)/\pi_1(M)$ where $\pi_1(M)$ acts on the first
factor as covering transformations and via the adjoint representation on the second factor.

We will be most interested in two cases: when $\rho$ is the holonomy representation of a hyperbolic structure on a 3-manifold $M$ or a projective structure on a surface $S$. In the former
case $E(M)$ is the bundle of germs of Killing fields on $M$ and for
the latter case $E(S)$ is the bundle of germs of projective vector
fields on $S$. 

The bundle $E(M)$ has a flat connection, $d$, which allows us to
define deRham cohomology groups. Our main goal of this paper is to
calculate $H^1(M;E(M))$ for a certain class of hyperbolic
structures. In \S\ref{representations} we will see that $H^1(M; E(M))$ is the tangent space of the space of hyperbolic cone-structures on $M$.

%
%

For the remainder of this section we restrict to the case where $E = E(M)$ is the bundle of a hyperbolic 3-manifold $M$.

A hyperbolic
structure on $M$ can be defined by a {\em developing map}, $\dev:\tilde{M}
\longrightarrow \hthree$, and a {\em holonomy representation},
$\rho:\pi_1(M) \longrightarrow 
\Isom^+\hthree$, where $\dev$ is a local homeomorphism that
commutes with the action of $\pi_1(M)$, i.e. $\dev(\gamma(x)) =
\rho(\gamma)\dev(x)$ for all $\gamma \in \pi_1(M)$ and $x \in
\tilde{M}$. 

A smooth 1-parameter family of hyperbolic structures, $M_t$, on
$M$ will have a smooth 1-parameter family of developing
maps, $\dev_t$. We call such a family a {\em local deformation} of the
hyperbolic structure. For each $x \in \tilde{M}$, $\dev_t(x)$ is a
smooth path in 
$\hthree$. The tangent vector at $\dev_0(x)$ can be pulled back by $\dev_0$
to a tangent vector at $x$ defining a vector field $v$ on
$\tilde{M}$. This vector field represents an {\em infinitesimal
deformation} of the hyperbolic structure.
Differentiating
$$\dev_t(\gamma(x)) = \rho_t(\gamma)\dev_t(x)$$
we see that 
\begin{equation}
\label{rhoderivative}
v - \gamma_*v = \dot{\rho}(\gamma)
\end{equation}
where $\dot{\rho}(\gamma)$ is the tangent vector of the path $\rho_t(\gamma)$ at $t=0$. The tangent vector of a path in $PSL_2\cx$ is an element the Lie algebra $sl_2\cx$ which is implicitly identified with a Killing field in \eqref{rhoderivative}. The Killing field $\dot{\rho}$ describes the {\em infinitesimal change in holonomy} of $\gamma$ induced by the deformation. In general if $v$ is a vector field on $\tilde{M}$ with $v - \gamma_* v$ a Killing field for all $\gamma \in \pi_1(M)$ then $v$ is {\em automorphic}.

A vector field on $M$ (or $\tilde{M}$) can be lifted to a section of $E$ (or $\tilde{E}$). In fact this can be done in three different ways which we now describe.
To do this we note that
the Lie algebra, $sl_2\cx$, has a complex structure that can be
geometrically interpreted using the $\curl$ operator on vector
fields. If $v \in sl_2\cx$ is a Killing vector field on $\hthree$ then
$\curl v$ will also be a Killing field and
$\curl \curl v = -v$. (The $\curl$ we are using differs from the usual curl be a factor if $-\frac{1}{2}$. We'll say more about this below.) We can then
define $\imath v = \curl v$. This will coincide with the usual complex
structure on $sl_2\cx$.

If $s$ is a section of $E$ and $p \in M$ then $s(p)$ is a Killing
field in a neighborhood of $p$ and $s(p)(p)$ will be a vector in the
tangent space of $M$ at $p$. We then define vector fields $\Re s$ and
$\Im s$ by $\Re s(p) = s(p)(p)$ and $\Im s(p) = -(\curl s(p))(p)$. Since
the Killing field $s(p)$ is uniquely determined by it value $s(p)(p)$ at $p$,
and the value $(\curl s(p))(p)$ of its $\curl$ at $p$ this
defines an isomorphism between $E$ and $TM \oplus TM$.

If $v$ is a vector field on $M$ we lift $v$ to sections $V$ and $\imath V$
of $E$ by setting $\Re V = v$ and $\Im V = 0$ while $\Re \imath V = 0$ and
$\Im \imath V = v$. Throughout the paper we will use this notational
convention of denoting vector fields on $M$ by lower case letters and
their corresponding sections of $E$ by uppercase letters.

We now define
one final method for lifting a vector field $v$ to a section of $E$.
For each point $p$ on $M$ we can find a Killing
field, $v_p$, in a neighborhood of $p$ that best approximates $v$ at
$p$. That is $v_p$ will be the unique Killing field such that $v_p(p) = v(p)$
and $(\curl v_p)(p) = (\curl v)(p)$. The {\em canonical lift} of $v$,
$s$, is defined by $s(p) = v_p$. Working through the definitions we see
that $(\curl s(p))(p) =
is(p)(p) = (\curl v)(p)$ so $\Re is = - \Im s = \curl v$ and $s=V -
\imath \curl V$. (In \cite{Hodgson:Kerckhoff:cone} this is the
definition of the canonical lift.) 

We define a section $s$ of $\tilde{E}$ to be {\em automorphic} if $s -
\gamma_* s$ is constant. An automorphic vector field and an
automorphic section both describe infinitesimal deformations of the
hyperbolic structure on $M$. We have the following relationship.
\begin{prop}
\label{canonicallift}
Let $s$ be the canonical lift of an automorphic vector field,
$v$. Then $s$ is an automorphic section.
\end{prop}

{\bf Proof.} By naturality $s - \gamma_* s$ will be
the canonical lift of $v- \gamma_* v$. Since $v$ is automorphic, $v -
\gamma_* v$ will be a Killing field. By definition the
canonical lift of a Killing field will be
constant. Therefore $s - \gamma_* s$ is constant and $s$ is
automorphic. \qed{canonicallift}

There are two simple but trivial ways to construct a local deformation of a hyperbolic structure. First we can post-compose the developing map with a smooth family $\sigma_t$ of isometries of $\hthree$ with $\sigma_0 = \id$. In this case the associated automorphic vector field will be a Killing field. The other method is to pre-compose the developing map with the lift of an isotopy of $M$. In this case the associated automorphic vector field will actually be an equivariant vector field. For this reason we say that an infinitesimal deformation is trivial if it is
the sum of a Killing field and an equivariant vector field. In terms
of sections a deformation is trivial if it is a constant section plus
an equivariant section. Two deformations are equivalent if they differ
by a trivial deformation. This definition holds for both vector fields
and sections.

If $s$ is an automorphic section then $\omega = ds$ will be an
equivariant 1-form because
$$ \omega - \gamma_*\omega =ds - \gamma_* ds = d(s-\gamma_*s) = 0 $$
since $s - \gamma_* s$ is constant.
Therefore $\omega$ descends to an $E$-valued 1-form on $M$. 
If $s$ is equivariant then $\omega$ will be an exact 1-form. If $s$ is
constant $\omega$ will be zero. Also $s_1$ and $s_2$ will be
equivalent deformations if and only if $\omega_1 = ds_1$ and $\omega_2
= ds_2$ differ by an exact 1-form. Therefore the deRham cohomology
group $H^1(M;E)$ describes the space of infinitesimal deformations.

{\bf Remark.} For a vector field $v$ on $\tilde{M}$ with canonical
lift $s$ the $E$-valued 1-form, $\omega = ds$ should be compared with
Thurston's description of the Schwarzian derivative of complex
analysis. In particular, if $f$ is a univalent holomorphic function
then for each $z$ we can find an {\em osculating M\"obius
transformation}, $M_z^f$, which is the unique M\"obius transformation
whose 2-jet agrees with the 2-jet of $f$ at $z$. The {\em Schwarzian
derivative} of $f$ is the derivative of $M_z^f$ in $PSL_2\cx$.  This
definition makes it apparent that the Schwarzian measures how far $f$
differs from a projective map just as $\omega$ measure how far $v$
differs from a Killing field.

\medskip
%
%
%
We now define a metric on $E$ and $\tilde{E}$. 
If $x \in \hthree$ and $v,w \in sl_2\cx$ we can
define an inner product on $sl_2\cx$ depending on $x$ by
$$ \langle v,w \rangle_x = \langle v(x), w(x) \rangle + \langle \imath
v(x), \imath w(x) \rangle $$ where $\langle,\rangle$ is the standard
inner product on $\hthree$. If $\gamma \in PSL_2\cx$ then 
\begin{equation}
\label{innerproduct}
\langle v,w \rangle_x = \langle \gamma_*v, \gamma_*w
\rangle_{\gamma(x)}
\end{equation}
 where $\gamma_*$ acts on $sl_2\cx$ by the adjoint
representation. Via the developing map this defines an inner product on the
fibers of $\tilde{E}$. By (\ref{innerproduct}) this inner product is
invariant under the action of $\pi_1(M)$ and therefore descends to an
inner product on the fibers of $E$.

The inner product determines a bundle map from $E$ to the dual
bundle $E^*$. If $\alpha$ is an $E$-valued form we write its
$E^*$-valued dual as $\alpha^\sharp$. For an $E^*$-valued form,
$\alpha$, the dual form is $\alpha^\flat$.
In local coordinates we can write any $k$-form as a sum of terms of the form $s \omega$ with $s$ an $E$-valued section and
$\omega$ a real valued $k$-form. We then use the Hodge $*$-operator for
the hyperbolic metric on real forms to define $*(s\omega) = s(*\omega)$
and $ * (s\omega)^\sharp =(s^\sharp)(*\omega)$ and extend the definition to a arbitrary $E$-valued $k$-form linearly. It is easy to see that this local
definition is well defined and this allows us to define an inner product
on $E$-valued $k$-forms $\alpha$ and $\beta$ by
\begin{equation}
\label{Emetric}
(\alpha, \beta) = \int_M \alpha \wedge * \beta^\sharp.
\end{equation}
Here the wedge product between an $E$-valued form and an $E^*$-valued
form is a real form. We also define
$$\|\alpha\|^2 = \alpha \wedge *\alpha^\sharp.$$

The bundle $E^*$ will also have a flat connection with exterior derivative $d^*$.
However, differentiating in $E^*$ is not the same as differentiating
in $E$. More explicitly let $\del \omega= (d*
\omega^\sharp)^\flat.$ We shall see 
shortly that $d \neq \del$. We use $\del$ to define a formal adjoint
for $d$. Let $\delta = (-1)^{k}*\del*$ where $\delta$ acts on
an $E$-valued $k$-form on an $3$-dimensional hyperbolic manifold. 
Then $(d\alpha, \beta)= (\alpha, \delta\beta)$ if $\alpha$ and
$\beta$ are $C^\infty$ $k$-forms with compact support, i.e. $\delta$
is the formal adjoint of $d$. We can now define the Laplacian,
$\Delta = d\delta + \delta d$.

In local coordinates there is a nice formula for $d$ and $\delta$ in terms of the Riemannian connection $\nabla$ and algebraic operators. If $\{e_i\}$ is an orthonormal frame field with dual co-frame field $\{\omega^i\}$ we have
\begin{equation}
\label{dform}
d= \sum_i \omega^i \wedge (\nabla_{e_i} + \ad(E_i))
\end{equation}
and
\begin{equation}
\label{deltaform}
\delta = - \sum_j i(e_j)(\nabla_{e_j} - \ad(E_j)).
\end{equation}
Here $i()$ is the interior product on forms. The operator $\ad(E_i)$ takes a Killing field $Y$ to the Killing field $[E_i, Y]$ where $[,]$ is the usual bracket on vector fields. (Recall that $E_i$ is the real lift of $e_i$.) We also need to decompose $d$ and $\delta$ into their real and symmetric parts. That is let $D = \Re d$, $T = \Im d$, $D^* = \Re \delta$ and $T^* = \Im \delta$. Note that $D^*$ and $T^*$ are the formal adjoints of $D$ and $T$ respectively. It is also worth noting that $\del = D - T$.

In \cite{Matsushima:Murakami:harmonic} it is calculated that $T^*D +
D^*T + TD^* + DT^* =0$ and therefore $\Delta$ is a real operator.
This leads to a Weitzenbock formula
\begin{equation}
\label{wietzenbockformula}
\Delta = \Delta_D + H
\end{equation}
where $\Delta_D = D^*D + DD^*$ and $H = T^*T +TT^*$. Note that $H$ is
a purely algebraic operator.



The tangent bundle, $TM$, has an exterior derivative $\hat{d}$ and the
hyperbolic metric gives a co-derivative, $\hat{\delta}$. For a vector
field $v$ let $\hat{v}$ be the dual 1-form. We use similar notation to define the real valued Laplacian. Namely let $\hat{\Delta} = \hat{d} \hat{\delta} + \hat{\delta} \hat{d}$. There is also a Weitzenb\"{o}ck formula relating the two Laplacians:
\begin{equation}
\label{delta=deltahat+4}
\widehat{\Delta V} = \hat{\Delta} \hat{v} + 4\hat{v}.
\end{equation}

The divergence and curl of $v$ can be defined in terms of $\hat{d}$ and $\hat{\delta}$. That is, $\div v =
\hat{\delta}\hat{v}$ is the {\em divergence} of $v$ and the {\em curl} of $v$ is defined by the formula $\widehat{\curl v} = -\frac{1}{2}*\hat{d}\hat{v}$. (One usually defines the curl to be the vector field dual to $*\hat{d}\hat{v}$. The factor of $-\frac{1}{2}$ is chosen such the curl of a Killing field is that same as multiplication by $\imath$ in $sl_2\cx$.) 

For a vector field $v$, $\nabla v$ is a
tensor of type (1,1), i.e. a section of the bundle $\Hom (TM,
TM)$. The divergence, curl and a third quantity, the strain, completely determine $\nabla v$. In particular, $\div v$ is the trace of $\nabla v$ and measures the infinitesimal change in volume. By definition the {\em strain}, $\str v$, of $v$ is the symmetric, traceless part of $\nabla v$. It measures the infinitesimal change in conformal structure. The divergence and strain
together measure the infinitesimal change in metric. The skew-symmetric
part of $\nabla v$, $\skew \nabla v$, is naturally identified with the curl of $v$. More explicitly there is an isomorphism from $TM$ to skew-symmetric sections of $\Hom(TM, TM)$.  We define this isomorphism by choosing an orthonormal from field $\{e_1, e_2, e_3\}$ with dual co-frame field $\{\omega^1, \omega^2, \omega^3\}$ and sending $e_i$ to $e_{i+1} \otimes \omega^{i+2} - e_{i+2} \otimes \omega^{i+1}$ (the indices are measured mod 3). Under this isomorphism $\skew \nabla v$ is exactly $\curl v$ as defined above.
The Riemannian
metric gives a norm to $\Hom(TM,TM)$ for which this decomposition is
orthogonal. 

%
%
%

The real and imaginary parts of an $E$-valued 1-form are both vector valued 1-forms or sections of $\Hom(TM, TM)$. Using the formulas above Hodgson and Kerckhoff relate the real and imaginary parts of an $E$-valued 1-form to the divergence, strain and curl of a vector field.  These results, which we summarize in the following theorem, can be found in \S2 of \cite{Hodgson:Kerckhoff:cone}. 

\begin{theorem}
\label{structure} 
Let $s$ be an automorphic section of $\tilde{E}$. Then there exists an automorphic vector field $v$ and an equivariant vector field $w$ such that $s = V - \imath \curl V + \imath W$. Moreover 
\begin{enumerate}
\item $\sym \Re ds = \nabla v$;

\item $\skew \Re ds = w$;

\item if $v$ is divergence free and harmonic and $w \equiv 0$ then $\Re ds = \str v$ and $\Im ds = -\str \curl v$.
\end{enumerate}
%
%
%
%
%
%
%
%
%
\end{theorem}

An $E$-valued 1-form $\omega$ is a {\em Hodge form} if there exists an automorphic, divergence free, harmonic vector field $v$ with canonical lift $s$ such that $\omega = ds$. There is a very simple formula for the $L^2$-norm of a Hodge form. It is essentially Proposition 1.3 of
\cite{Hodgson:Kerckhoff:cone}.
\begin{prop}
\label{boundaryterm}
Let $M$ be a compact hyperbolic manifold with boundary and $\omega$ a Hodge form on $M$. Then
$$\int_M \|\omega\|^2 = \frac{1}{2} \int_{\del M} \imath \omega \wedge \omega^\sharp$$
where $\del M$ is oriented with inward pointing normal.
\end{prop}

\section{Geometrically finite ends}
\label{ends}

\subsection{Projective structures and geometrically finite ends}

Throughout this section we let $M = S \times [0,\infty)$ where $S$ is a
closed surface of genus $>1$. We also assume that $M$ has a complete
hyperbolic structure with boundary $S \times \{0\}$.

A {\em projective structure}, $\Sigma$, on a surface is
given by an atlas of charts with 
image in $\chat$ and transition maps M\"obius or projective
transformations. As with hyperbolic structures, a projective structure
can be given by a developing map and a holonomy representation.
If $\Sigma_1$ and $\Sigma_2$ are projective structures then
$\Sigma_1 \cong \Sigma_2$ if there exists a projective homeomorphism
from $\Sigma_1$ to $\Sigma_2$.

\begin{defn}
\label{finiteends}
$M$ is a {\em geometrically finite end without rank one cusps} if it is compactified by a projective structure $\Sigma$ on $S \times \{\infty\}$. Then $\Sigma$ is the {\em projective boundary} of $M$.
\end{defn}

Since we will not discuss rank one cusps in this paper we will simply refer to such ends as geometrically finite.

To see this definition more explicitly we recall that $\chat$ naturally compactifies $\hthree$. We refer to this compactification as $\hbar = \hthree \cup \chat$. Then $PSL_2\cx$ acts continuously on $\hbar$ as isometries of $\hthree$ and projective transformations of $\chat$. Then $M$ is geometrically finite if it has an atlas of hyperbolic charts that extends continuously to a atlas for a projective structure on $S \times \{\infty\}$. In fact if $M$ is geometrically finite there will be a developing map
$$D: \tilde{S} \times [0, \infty] \longrightarrow \hbar$$
that restricts to a developing map for the hyperbolic structure $M$ on $\tilde{S} \times [0, \infty)$ and a developing map for the projective structure $\Sigma$ on $\tilde{S} \times \{\infty\}$.

We refer
to the bundles $E(\Sigma)$ and $E(\tilde{\Sigma})$ as $E_\infty$ and
$\tilde{E}_\infty$, respectively. $E_\infty$ and $\tilde{E}_\infty$
are the bundles of germs of projective vector fields on $\Sigma$ and
$\tilde{\Sigma}$ respectively.

The product structure on $M \cup \Sigma$ allows us to define a projection map
$$\Pi: M \longrightarrow \Sigma$$
by the formula $\Pi(p,t) = p$. Let $\tilde{\Pi}$ be a lift of $\Pi$ to the universal covers $\tilde{M}$ and $\tilde{\Sigma}$. The corresponding bundles $\tilde{E}$ and $\tilde{E}_\infty$ have canonical product structures so we can use $\tilde{\Pi}$ to pull back sections of $\tilde{E}_\infty$ to sections of $\tilde{E}$. If we restrict $\tilde{\Pi}$ to $\tilde{S} \times \{0\}$ in $\tilde{M}$ we can also push forward sections of $\tilde{E}$ to sections of $\tilde{E}_\infty$.

\begin{lemma}
\label{preserve automorphic} Let $s$ be an automorphic section of
$\tilde{E}$ and $s_\infty$ an automorphic section of
$\tilde{E}_\infty$. 
\begin{enumerate}
\item $\tilde{\Pi}_* s$ and $\tilde{\Pi}^* s_\infty$ will be
equivariant iff $s$ and $s_\infty$, respectively, are equivariant.

\item $\tilde{\Pi}_* s$ and $\tilde{\Pi}^* s_\infty$ are
automorphic. 

\item Automorphic sections $s'$ and
$s_\infty'$  will be equivalent as infinitesimal deformations to $s$ and
$s_\infty$, respectively, iff $\tilde{\Pi}_* s'$ and $\tilde{\Pi}^*
s_\infty'$ are 
equivalent to $\tilde{\Pi}_* s$ and $\tilde{\Pi}^* s_\infty$, respectively.
\end{enumerate}

Therefore there are isomorphisms $\Pi_*:H^1(M; E) \longrightarrow
H^1(S; E_\infty)$  and $\Pi^*:H^1(S; E_\infty) \longrightarrow
H^1(M; E)$.

\end{lemma}

{\bf Proof.} The actions of $\pi_1(M) = \pi_1(S)$ on $\tilde{E}$ and
$\tilde{E}_\infty$, respectively, will commute with $\Pi$ so 
$$\tilde{\Pi}_* s - \gamma_*\tilde{\Pi}_*s = \tilde{\Pi}_*(s -
\gamma_* s) $$
for all $\gamma \in \pi_1(M)$. This implies (1) for $\tilde{\Pi}_*$ and
(1) then implies (2) and (3). The proof for $\tilde{\Pi}^*$ is similar.

Given a closed $E$-valued 1-form $\omega$ we can integrate $\omega$ to
find an automorphic section, $s$, of $\tilde{E}$ such that $ds =
\omega$. We then define $\Pi_* \omega = d \Pi_* s$.  We similarly
define $\Pi^*$ for and $E_\infty$-valued 1-form. (1), (2) and (3)
imply that this defines maps between $H^1(M;E)$ and
$H^1(S;E_\infty)$. If $s' = \Pi^*\Pi_* s$ then $s=s'$ on $\Sigma_0$ so
$ds'$ is cohomologous to $ds$. Therefore the map $\Pi_*$ is one-to-one
from $H^1(M;E)$ to $H^1(S;E_\infty)$ and by similar reasoning $\Pi^*$
is one-to-one from $H^1(S;E_\infty)$ to $H^1(M;E)$ and both
$\Pi^*\Pi_*$ and $\Pi_*\Pi^*$ are the identity map. Therefore the maps
$\Pi_*$ and $\Pi^*$ are isomorphisms. \qed{preserve automorphic}

\medskip

{\bf Remark.} In general a vector field on a
geometrically finite end will not extend continuously to the conformal
boundary. What Lemma \ref{preserve automorphic} allows us to do is
replace on automorphic vector, $v$, on $\tilde{M}$ with an equivalent
vector field that does extend continuously. Namely if $s$ is the
canonical lift of $v$ then $v'=\Re\Pi^*\Pi_*s$ will be equivalent to
$v$ and $v' \cup \Re \Pi_* s$ will be continuous on $\tilde{M} \cup
\tilde{S}$.

%
\subsection{Extending sections via horosphere projections}
\label{vector fields on ends}

%
%
We will need to extend automorphic vector fields on $\tilde{\Sigma}$ to harmonic, divergence free vector fields on $\tilde{M}$. We first describe a method for extending a vector field on $\cx$ to $\hthree$. More precisely we extend a vector field to a harmonic section of $E(\hthree)$.

We will use the following orthonormal frame field on $\hthree$:
Working in the upper half space model, $\hthree = \{(x,y,t):
t >0 \}$, let $e_1 = t \frac{\del}{\del x}$, $e_2 = t
\frac{\del}{\del y}$ and $e_3 = t \frac{\del}{\del t}$
with corresponding real $E$-valued sections $E_i$.
Also let $R_i = \imath E_i$ and
let $\omega^i$ be the dual real 1-form for $e_i$. Using (\ref{dform}) and (\ref{deltaform}) we see that
\begin{eqnarray}
dE_1 & = & E_3\omega^1 + R_3\omega^2 - R_2\omega^3 \nonumber \\
dE_2 & = & -R_3\omega^1 + E_3\omega^2 + R_1\omega^3 \label{dcalc}\\
dE_3 & = & -(E_1 - R_2)\omega^1 - (R_1+E_2)\omega^2 \nonumber
\end{eqnarray}
and
\begin{eqnarray}
\partial E_1 & = & E_3\omega^1 - R_3\omega^2 + R_2\omega^3 \nonumber \\
\partial E_2 & = & R_3\omega^1 + E_3\omega^2 - R_1\omega^3
\label{delcalc} \\ 
\partial E_3 & = & -(E_1 + R_2)\omega^1 + (R_1-E_2)\omega^2. \nonumber
\end{eqnarray}

If $v$ is a projective vector field on $\cx \subset \chat$ then there is an obvious way to extend $v$ to a harmonic vector field on $\hthree$. Namely there is a unique Killing field on $\hthree$ that extends continuously to $v$ on $\cx$. Every projective vector field is of the form $v(z) = p(z) \ddz$ where $p$ is a quadratic polynomial. Then the section
\begin{equation}
\label{extkilling}
\frac{p(w)}{t}(E_1 - R_2) + p_z(w) E_3 - \frac{t p_{zz}(w)}{2}(E_1 + R_2)
\end{equation}
is constant and evaluates at every point to the Killing field that extends continuously to $v$.

Now let $v(z) = f(z) \ddz$ be an arbitrary smooth vector field on $\cx$. Then the {\em canonical lift}, $s_\infty(z)$, of $v$ is the section of $E(\cx)$ defined by the formula
\begin{equation}
\label{liftonchat}
s_\infty(w) = \left(f(w) + f_z(w)(z-w) + \frac{f_{zz}(w)}{2}(z-w)^2
\right)\ddz.
\end{equation}
Similar to the canonical lift of vector fields on hyperbolic space, $s_\infty(w)$ is the projective vector field that best approximates $v$ at $w$.
Next we define a section of $E(\hthree)$ by the formula $s(w,t)  = s_\infty(w)$. Using (\ref{extkilling}) we see that
\begin{equation}
\label{horosphere extension}
s(w,t)  = \frac{f(w)}{t}(E_1 - R_2) + f_z(w)E_3 -
\frac{tf_{zz}(w)}{2}(E_1 + R_2).
\end{equation}
Note that the vector fields $\Re s$ and $\Im s$ extend continuously to $v$ and $- \imath v$, respectively.
%
%

To calculate $\Delta s$ we use the formula,
\begin{equation}
\label{product formula}
\Delta (fs) = (\hat{\Delta} f) s - 2*(*\hat{d}f \wedge Ds) +f(\Delta s).
\end{equation}
where $f$ is a complex valued function and $s$ an $E$-valued section.
If $f$ is a function in the upper half space model of $\hthree$ then
\begin{equation}
\label{df}
\hat{d}f = tf_x \omega^1 + tf_y \omega^2 + tf_t \omega^3 
\end{equation}
and
\begin{equation}
\label{deltaf}
\hat{\Delta} f = tf_t - t^2(f_{xx} + f_{yy} + f_{tt}).
\end{equation}

After a straightforward but long calculation, using these formulas, we
have: 
\begin{equation}
\label{laplacian of horosphere extension}
\Delta s(w,t) = -2tf_{z \zbar}(w)(E_1-R_2) - 2t^2f_{zz \zbar}(w)E_3 +
2t^3f_{zzz\zbar}(w)(E_1+R_2).
\end{equation}

There are several things to notice in this formula. If $v = f \ddz$ is
conformal, i.e. f is a holomorphic function, then $\Delta s = 0$. If
$v$ is not conformal then $\|\Delta s \| \rightarrow 0$ as $t
\rightarrow 0$. In other words as we approach the ideal boundary the
norm of $\Delta s$ goes to zero. In fact we can estimate the rate of
decrease. Let $d(w,t) = -\log t$ be the distance of a point from the
horosphere $t=1$. Then $\| \Delta s(w,t) \|  = 4|f_{z
\zbar}(z)|e^{-d(w,t)} + o(t^2)$. ($o(t)$ is a function such that
$\frac{o(t)}{t}$ is bounded as $t \rightarrow 0$.)

For a conformal vector field there is also a nice expression for $ds$.
Namely
\begin{equation}
\label{ds for horosphere extension}
ds = -\frac{t^2 f_{zzz}}{2}(E_1 + R_2)(\omega^1 + \imath\omega^2).
\end{equation}
Here we note that $\Re ds$ and $\Im ds$ are symmetric and traceless,
therefore they are the strains of the divergence free vector fields
$\Re s$ and $\Im s$, respectively. We can also calculate the norm,
$\|ds(p,t)\| = |f_{zzz}(w)|e^{-2d}$, which should be thought of as the
norm of the strains of $\Re s$ and $\Im s$.

{\bf Remark.} By differentiating (\ref{liftonchat}) we see that
$$ds_\infty(w) = \left(\frac{f_{zzz}(w)}{2}(z-w)^2\ddz\right)dw.$$
The quantity $f_{zzz}(w)$ is the infinitesimal version of the
Schwarzian derivative discussed in the remark in \S\ref{background}. In particular we see that the Schwarzian derivative of the
vector field at infinity determines the strain of the extended vector
field in hyperbolic space.

\subsection{Convex parallel surfaces in $\hthree$}

We now describe some results on parallel surfaces in hyperbolic
space. A detailed study of such surfaces can be found in
\cite{Epstein:horospheres}. All of the results described in this
section can be found in \S3 of this paper.

Let $S$ be a smooth convex surface in $\hthree$. For a point $p \in S$ normalize the position of $S$ in $\hthree$ such that $p = (0,0,1)$ in the upper half space model and such that the directions of principal curvature at $p$ are $(1,0,0)$ and $(0,1,0)$. In a neighborhood $V$ of the geodesic ray normal to $S$ at $p$ we can foliate $\hthree$ by surfaces $S_t$ such that each $S_t$ is  equidistant from $S$ and contains the point $(0,0,t)$. Note that the distance between $S$ and $S_t$ is $- \log t$. There are maps $\pi_t : S \longrightarrow S_t$ that take each point in $S$ to the nearest point in $S_t$. There 
is also a projection $\Pi: V \longrightarrow \cx$ that takes each point on $S_t$ to the limit of its normal ray in $\cx$. Note that $\Pi$ restricted to $S$ is the limit of the maps $\pi_t$ as $t \rightarrow 0$.

Let $k_1(t)$ and $k_2(t)$ be the principal curvatures of the surface $S_t$ at the point $(0,t)$. Then
\begin{equation}
\label{curvature calc}
k_i(t) = \frac{1 + k_i(0) + t^2(k_i(0) - 1)}{1 + k_i(0) + t^2(1 - k_i(0))}.
\end{equation}

We also need to calculate the derivatives of the $\pi_t$ and $\Pi$. We will express the answer
in terms of the orthonormal frame field $\{e_1, e_2, e_3\}$ for $\hthree$ defined in the previous section and the frame field $\left\{ \frac{\del}{\del x} , \frac{\del}{\del y} \right\}$ for $\cx$. Then the derivative of $\pi_t$ at $(0,0,1)$ is 
\begin{equation}
(\pi_t)_* = \left( \begin{array}{cc}
\frac{1 + k_1(0) + t^2(1-k_1(0))}{2t} & 0 \label{vperp}\\
0 & \frac{1+k_2(0) + t^2(1- k_2(0))}{2t}
\end{array} \right)
\end{equation} 
and the derivative of $\Pi$ at $(0,0,t)$ is
\begin{equation}
\label{derivatives}
\Pi_* = \left( \begin{array}{ccc}
\frac{(1+k_1(t))t}{2} & 0 & 0 \\
0 & \frac{(1+k_2(t))t}{2} & 0 
\end{array} \right).
\end{equation}

%

%
\subsection{Extending sections via convex surfaces}
We continue to assume that $S$ be an embedded convex surface in $\hthree$ with position normalized as in the previous section. However, in this section it will be convenient to use complex coordinates for the upper half space model. That is we view the first two real coordinates as a single complex coordinate. 

Let $v$ be a conformal vector field on $\cx$ and let $s_ \infty$ be the canonical lift of $v$. Define a section of $E(\hthree)$ by the formula $s = \Pi^* s_\infty$. We will calculate $ds$ and $\Delta s$ along the ray $(0,t)$.

\begin{prop}
\label{decayrate}
$\|ds(0,t)\|$, $\|(\Delta s)(0,t)\|$, $\|(\div \Re s)(0,t)\|$ and
$\|\hat{d}(\div \Re s)(0,t)\| = o(t^2)$.
\end{prop}

{\bf Proof.} For the special case when $S$ is a horosphere we calculated $s$, $ds$ and $\Delta s$ in \S\ref{vector fields on ends}. In general $s$ will be the sum of the horosphere extension, $s_h$, and a correction term $s_c$.

Let $z = \Pi(w,t)$. Using (\ref{extkilling}) we see
\begin{eqnarray}
s(w,t) & = & \left[f(z) + f_z(z)(w-z) +
\frac{f_{zz}(z)(w-z)^2}{2} \right]\frac{E_1 - R_2}{t} \nonumber \\
& & + \left[f_z(z) + f_{zz}(z)(w-z) \right]E_3 -
f_{zz}(z)\frac{t(E_1 + R_2)}{2} \nonumber \\
& = & [f(w) - g_1(w,z)]\frac{E_1-R_2}{t}  \nonumber \\
& & + [f_z(w) - g_2(w,z)]E_3 \nonumber \\
& & - [f_{zz}(w) - g_3(w,z)]\frac{t(E_1+R_2)}{2} \nonumber
\end{eqnarray}
where
\begin{eqnarray}
\label{generals}
g_1(w,z) & = & (w-z)^3 \sum_{n=3}^\infty
\frac{f^{(n)}(z)}{n!}(w-z)^{n-3} \nonumber \\
g_2(w,z) & = & (w-z)^2 \sum_{n=2}^\infty
\frac{f^{(n+1)}(z)}{n!}(w-z)^{n-2} \\
g_3(w,z) & = & (w-z) \sum_{n=1}^\infty
\frac{f^{(n+2)}(z)}{n!}(w-z)^{n-1}. \nonumber
\end{eqnarray}
Therefore $s = s_h - s_c$ where,
\begin{equation}
s_h(w,t) = \frac{f(w)}{t}(E_1 - R_2) + f_z(w) E_3 -
\frac{tf_{zz}(w)}{2}(E_1 + R_2),
\end{equation}
and,
\begin{equation}
\label{correct}
s_c(w,t) = \frac{g_1(w,z)}{t}(E_1 - R_2) + g_2(w,z)E_3 -
\frac{tg_3(w,z)}{2}(E_1 + R_2).
\end{equation}
We have already calculated $ds_h$ and $\Delta s_h$ so we are left to
calculate $ds_c$ and $\Delta s_c$.

Let $G_i(w,t) = g_i(w, \Pi(w,t))$ for $i=1,2,3$.
By our normalization $\Pi(0,t) = 0$ for all $t$,
so $w-z=0$ when $w =0 $. Furthermore the $G_i$ extend to smoothly to
$\cx$ with $G_i(w,0) = g_i(w, \Pi(w,0)) = g_i(w,w) = 0$. For this
reason the Euclidean derivatives of the $G_i$ will be bounded on
$(0,t)$. Then (\ref{df}) and (\ref{deltaf}) imply that
$\|\hat{d}G_i(0,t)\| = o(t)$ and $\|\hat{\Delta} G_i(0,t)\| = o(t^2)$.

In fact using the product formula for the real Laplacian,
\begin{equation}
\label{product formula for functions}
\hat{\Delta} (fg) = (\hat{\Delta} f) g - 2<\hat{d}f,\hat{d}g> + f
(\hat{\Delta} g),
\end{equation}
and the Leibniz rule we obtain
\begin{equation}
\begin{array}{cclclcc}
G_1(0,t) &=& \hat{d}G_1(0,t) &=& \hat{\Delta} G_1(0,t) &=& 0  \\
G_2(0,t) &=& \hat{d}G_2(0,t)&=&0& &  \label{leibniz} \\
G_3(0,t) &=&0.& & & &  
\end{array}
\end{equation}

(\ref{leibniz}) implies that
\begin{equation}
\label{dsc}
ds_c(0,t) = -\frac{t}{2}(E_1 + R_2)\hat{d}G_3(0,t)
\end{equation}
so $\|ds_c(0,t) \| = o(t^2)$. We saw in (\ref{ds for horosphere
extension}) that $\|ds_h(0,t)\| = o(t^2)$ so $\|ds(0,t)\| =
o(t^2)$. Furthermore $\div \Re s$ is the trace of $\Re ds$, so $\|(\div
\Re s)(0,t)\| = o(t^2)$.

To estimate $\hat{d}(\div \Re s)$ we will need to know $\div \Re s$
more explicitly. Here we use (\ref{derivatives}) from the previous
section to find
\begin{equation}
\label{dg3}
\hat{d}G_3(0,t) = \frac{tf_{zzz}(0)}{2}\left((1-k_1(t))\omega^1 +
\imath (1-k_2(t))\omega^2 \right).
\end{equation}
Therefore
$$ds_c(0,t) = -\frac{t^2f_{zzz}(0)}{4}(E_1 + R_2)\left((1-k_1(t))\omega^1 + \imath (1-k_2(t))\omega^2 \right)$$
and
\begin{equation}
\label{div}
(\div \Re s_c)(0,t) = \frac{t^2 \Re f_{zzz}(0)}{4}(k_1(t)-k_2(t)).
\end{equation}
Since $\div \Re s_h = 0$, $\div \Re s = \div \Re s_c$ and therefore by
(\ref{div})
$(\div \Re s)/t$ will extend to a smooth function on $\cx$ in a
neighborhood of zero. We again apply (\ref{df}) to see that
$\|\hat{d}(\div \Re s)(0,t)\| = o(t^2)$.

We are now left to calculate $(\Delta s_c)(0,t)$. By (\ref{product
formula}) and (\ref{leibniz})
\begin{eqnarray}
(\Delta s_c)(0,t) & = & (\hat{\Delta} G_2(0,t))E_3 + (\hat{\Delta}
G_3(0,t))\frac{t(E_1 + R_2)}{2} \nonumber \\
& &- 2*\left(*\hat{d}G_3(0,t) \wedge
D\left(\frac{t(E_1+R_2)}{2}\right)\right). 
\end{eqnarray}
We have estimated every term on the right except
$D\left(\frac{t(E_1+R_2)}{2}\right)$ which we can calculate using \eqref{dcalc} to see that it is $o(t)$. Therefore $\|(\Delta s_c)(0,t)\| = o(t^2)$ and
since $\Delta s_h = 0$, $\|(\Delta s)(0,t)\| =
o(t^2)$. \qed{decayrate}

\subsection{Extending vector fields on the boundary of a geometrically
finite end}

We now return to our geometrically finite end,
$M$. To apply the results of the previous two sections we need the following theorem:
\begin{theorem}
\label{convex}
If $M$ is a geometrically finite end then there
is a smooth, embedded, convex surface in $M$ whose inclusion
is a homotopy equivalence.
\end{theorem}

{\bf Proof.} Let $\Sigma$ be the projective boundary of $M$. For every projective structure there is a conformal developing map $f: U \longrightarrow \chat$ where $U$ is the upper half plane in $\cx$. We will extend $f$ to a developing map for $M$.

Let $M^f_z$ be the unique M\"obius transformation whose 2-jet agrees with $f$ at $z$. Let $P$ be the hyperbolic plane whose boundary is the real line on $\chat$. For each $p \in \hthree$ there is a unique geodesic $g$ through $p$ which is orthogonal to $P$. The geodesic $g$ will have exactly one endpoint $z \in U$. Define a map $F: \hthree \longrightarrow \hthree$ by the formula $F(p) = M^f_z(p)$. The construction is natural so $F$ will be a developing map for a geometrically finite end if $F$ is a local diffeomorphism. Furthermore $F$ extends continuously to $f$ so this end will have projective boundary $\Sigma$.

The upper half plane $U$ and the hyperbolic plane $P$ bound a half space $H$ in $\hthree$. Let $H_d$ be the set of points in $H$ whose distance from $P$ is $\geq d$. In \S3 of \cite{Anderson:projective} the derivative of $F$ is calculated. In particular it is shown that for $d$ sufficiently large $F$ is a local diffeomorphism restricted to $H_d$ (see p. 35). By Proposition 3.13 and Theorem 3.17, $F(\del H_d)$ will be convex for $d$ sufficiently large. For such $d$, $F$ restricted to $H_d$ will be a developing map for a geometrically finite end $M_d$ with projective boundary $\Sigma$. The boundary of $M_d$ will be a convex surface.

The isomorphism from the projective boundary of $M_d$ and the projective boundary of $M$ will extend to an isometry from $M_d$ to $M$. Although this isometry may not be defined on all of $M_d$ it can be defined on all but a compact submanifold of $M_d$. Therefore we can choose $d$ even larger so that this isometry is defined on all $M_d$. The image of $\del M_d$ under this isometry will be the desired convex surface. \qed{convex}

The convex surface we have constructed separates $M$ into a compact piece and a non-compact piece. The outward normal of the convex surface points into the non-compact piece so if we remove the compact piece $M$ will have concave boundary.
From now on we will assume this is the case. That is $\del M = S \times \{0\}$ is concave. We also assume that the product structure is chosen such that for a fixed $p \in S$ the path $\{p\} \times [0,\infty)$ is a geodesic ray in $M$ normal to $S \times \{0\}$. Furthermore we assume that the second parameter is a unit speed parameterization of the geodesic. This implies that $S \times \{t\}$ is convex for all $t >0$.

If $v$ is an automorphic vector field on $\tilde{\Sigma}$ we can
define a canonical lift of $v$ on a local projective chart as in
(\ref{liftonchat}). One can
then check to see that this lift is independent of the choice of chart
so $v$ has a well defined canonical lift, $s_\infty$, to all of
$\tilde{S}$. Furthermore if $v$ is automorphic $s_\infty$ will also be automorphic. Let $s = \Pi^* s_\infty$. Then Lemma \ref{preserve automorphic} implies that $s$ is also automorphic and therefore $ds$, $\Delta s$ and $\div \Re s$ will be equivariant and descend to objects of the appropriate type on $M$.

\begin{theorem}
\label{all l2} If $v$ is a conformal vector field then
$ds$, $\Delta s$, $\div \Re s$, and $\hat{d}(\div \Re s)$  all have
finite $L^2$ norm on $M$.
\end{theorem}

{\bf Proof.} Let $p \in M$ lie on $S \times \{t\}$. Proposition \ref{decayrate} implies that there exists
a continuous function $K: S \longrightarrow \reals^+$ such that $\|
ds(p) \| < K(\Pi(p))e^{-2t}$. Since $S$ is compact $K$ will bounded by some $K_1 >0$. 

Let $dA_t$ be the area form for $S \times \{t\}$. Note
that by (\ref{curvature calc}) and (\ref{vperp}), $\area (S \times \{t\}) < K_2
e^{2t} \area (S \times \{0\})$ where $K_2$ is determined by the maximal principal
curvature on $S \times \{0\}$. Then
\begin{eqnarray}
\int_M \|ds(p)\|^2
& = & \int_0^\infty \int_{\Sigma_t} \|ds(p)\|^2 dA_t dt \nonumber \\
& < & \int_0^\infty (K_1)^2e^{-4t} K_2 \area (S \times \{0\})e^{2t} dt \nonumber \\
& = & \int_0^\infty (K_1)^2K_2 \area(S \times \{0\})e^{-2t} dt < \infty \nonumber
\end{eqnarray}
and $ds$ is in $L^2$. 

The proof for the other terms is similar. \qed{all l2}

\subsection{Harmonic deformations of rank two cusps}
\label{cuspdef}

A rank two cusp is the quotient of a horoball by a $\integers \oplus \integers$ group of parabolic isometries. Again it is convenient to work in the upper half space model where isometries can be defined by there action on $\chat$. Since any two parabolics of $\hthree$ are conjugate we can explicitly describe any cusp as the quotient of a horoball based at infinity by parabolics $\gamma_1(z) = z + 1$ and $\gamma_2(z) = z + \tau$ with $\Im \tau >0$. The cusp $M$ is homeomorphic to $T \times [0, \infty)$ where $T$ is a torus. We can choose this product structure such that each $T \times \{t\}$ is the quotient of a horosphere and therefore has an induced Euclidean metric. The conformal class of these metrics will be constant for $T \times \{t\}$ and it is determined by the Teichm\"uller parameter $\tau$.

We defined $\gamma_1$ and $\gamma_2$ by their action on $\cx$. The quotient of this action is a projective structure $\Sigma$ on the torus. To construct models for deformations of the rank two cusp we will first describe deformations of $\Sigma$.

Since the Euler characteristic of $\Sigma$ is 0 Poincare duality implies that the complex dimension of $H^1(\Sigma; E(\Sigma))$ is 2. An automorphic vector field on $\tilde{\Sigma}$ determines a cohomology class in $H^1(\Sigma; E(\Sigma))$. We claim that the automorphic vector fields $v_1 = \frac{z-\zbar}{2} \ddz$ and $v_2 = \frac{z^3 - z}{6} \ddz$ determine cohomology classes that are a basis for $H^1(\Sigma; E(\Sigma))$.

We first examine $v_1$. Note that $\tilde{\Sigma} = \cx$ where $\cx$ has the natural projective structure it inherits as a subset of $\chat$. It is then easy to check that $v_1$ is automorphic. In particular $v_1 - (\gamma_1)_* v_1 =  0$ and $v_1 - (\gamma_2)_* v_1 = \Im \tau \ddz$ are projective vector fields. If the cohomology class $v_1$ generates is trivial then $v_1$ is the sum of an equivariant vector field and a projective vector field. For this to be true there must be a projective vector field $v$ with $v(0) = v(1) = 0$ and $v(\tau) = v(\tau + 1) =  \Im \tau$. Since no such $v$ exist the cohomology class is non-trivial.

Similar reasoning applied to $v_2$ and any non-zero linear combination of $v_1$ and $v_2$ shows that $v_2$ also generates a non-trivial cohomology class and together $v_1$ and $v_2$ determine a basis of $H^1(\Sigma; E(\Sigma))$.

Just as in \S\ref{vector fields on ends} we can extend $v_1$ and $v_2$ to sections of $E(\hthree)$. By (\ref{horosphere extension}), $v_1$ extends to
$$s_1(w,t) = \frac{w - \bar{w}}{2t}(E_1 - R_2) + \frac{1}{2}E_3$$
and $v_2$ extends to
$$s_2(w,t) = \frac{w^3 - w}{6t}(E_1 - R_2) + \frac{3w^2 -1}{6}E_3 - \frac{wt}{2}(E_1 + R_2).$$
Both $s_1$ and $s_2$ will be automorphic sections with respect to the action of $\gamma_1$ and $\gamma_2$ on $\hthree$. Therefore
$$ds_1 = -\frac{1}{2}(E_1 - R_2)(\omega^1 - \imath \omega^2)$$
and
$$ds_2 = -\frac{t^2}{2}(E_1 + R_2)(\omega^1 + \imath \omega^2)$$
are equivariant and restrict to $E$-valued 1-forms on $M$. By Proposition \ref{preserve automorphic} there is an isomorphism between $H^1(\Sigma; E(\Sigma))$ and $H^1(M; E)$ and therefore $ds_1$ and $ds_2$ are a basis for $H^1(M; E)$.

We have shown that every cohomology class in $H^1(M;E)$ has a representative of the form
$$\omega = -\frac{b_1}{2}(E_1 - R_2)(\omega^1 - \imath \omega^2) - \frac{b_2 t^2}{2}(E_1 + R_2)(\omega^1 + \imath \omega^2).$$

\begin{prop}
\label{cusprep} The $E$-valued 1-form, $\omega$, is closed, co-closed
and traceless. Furthermore $\omega$ is in $L^2$ if and only if $b_2 = 0$ and if $\omega$ is in $L^2$ then
$$\int_M \|\omega\|^2 = \frac{|b_1|^2}{2} \operatorname{Area}(\del M).$$
\end{prop}

{\bf Proof.} By construction $\omega$ is closed. From our explicit description of $\omega$ we see that it is traceless. By (\ref{laplacian of horosphere extension}) we see that $\Delta s_1 =  \delta d s_1= 0$ and $\Delta s_2 = \delta ds_2 =0$. Therefore $\omega$ is co-closed.

The pointwise norm of $\omega$ is $\|\omega(w,t)\|^2 = |b_1|^2 + t^4 |b_2|^2$. The last two facts follow from integrating this norm over $M$.
\qed{cusprep}

We next describe the infinitesimal change in holonomy determined by $\omega$. Again it is easier to work with the projective structure on the torus and then use the isomorphism between $H^1(\Sigma; E(\Sigma))$ and $H^1(M; E)$. Let $\gamma_t$ be a smooth path in $PSL_2\cx$ with $\gamma_0(z)  =z + \beta$. Then the derivative of $\gamma_t$ at $t=0$ will be a projective vector field $(a_0 + a_1 z + a_2 z^2)\ddz$.  A straightforward calculation shows that the derivative of the trace of $\gamma_t$ at $t=0$ is $-\beta a_2$.
Therefore if $v$ is a vector field on $\cx$ automorphic with respect to $\gamma_0$ then the infinitesimal change in trace is determined by the $z^2$-coefficient of the projective vector field   $v - (\gamma_0)_* v = (a_0 + a_1 z + a_2 z^2) \ddz$. 

We now apply this to the vector  fields $v_1$ and $v_2$. For both $\gamma_1$ and $\gamma_2$ the $z^2$-coefficient of $v_1 - (\gamma_i)_* v_1$ is zero. Therefore the infinitesimal change in trace determined by $v_1$ is zero. For $v_2$ the $z^2$-coefficient of $v_2 - (\gamma_1)_*v_2$ is $\frac{1}{2}$ so the infinitesimal change in trace of $\gamma_1$ is $-\frac{1}{2}$. The $z^2$-coefficient of $v_2 - (\gamma_2)_* v_2$ is $\frac{\tau}{2}$ so the infinitesimal change in trace of $\gamma_2$ is $-\frac{\tau^2}{2}$.

Although $v_1$ does not change the holonomy of either $\gamma_1$ or $\gamma_2$ there is an infinitesimal change in the projective structure. In particular there is an infinitesimal change in the conformal structure. Recall that the Teichm\"uller space of the torus can be identified with the upper half plane $U = \{z \in \cx | \Im z > 0\}$. Any point $\tau \in U$ determines a parallelogram with vertices $0$, $1$, $\tau$ and $\tau + 1$. By identifying opposite sides of this parallelogram we obtain a conformal structure on the torus. The affine vector field $v_1$ fixes the side between $0$ and $1$ so the infinitesimal change in the Teichm\"uller parameter is given by $v_1(\tau) = \Im \tau$. Note that the Teichm\"uller metric on the Teichm\"uller space of the torus is the hyperbolic metric. In this metric the length of the vector $v_1(\tau)$ is $1$. In particular its length does not depend on $\tau$. 

\section{Hodge theory of deformations}
\label{hodge theory}

We are now ready to begin our analysis of geometrically finite hyperbolic cone-manifolds. We begin with some definitions. Let $N$ be a compact 3-manifold with boundary and let $\cC$ be a collection of simple closed curves in the interior of $N$. Let $M$ be the interior of $N - \cC$. A  singular metric $g$ on $\operatorname{int}N$ is a {\em hyperbolic cone-metric} if $g$ is smooth metric of constant sectional curvature $\equiv -1$ on $M$ while in neighborhood of a point $p \in \cC$ the metric has the form
$$dr^2 + \sinh^2 rd\theta^2 + \cosh^2 r dz^2$$
with $\theta$ measured modulo some $\alpha >0$. On each component of $c$ of $\cC$, $\alpha$ will be constant. Then $\alpha$ is the cone angle of the cone singularity at $c$. We further say that $g$ is geometrically finite (without rank one cusps) if $g$ extends to a projective structure on each non-toral component of $\del N$.

The complement of any compact core of $M$ will contain ends of three types: geometrically finite ends, rank two cusps and neighborhoods of the cone singularity.
%
%
%
%
For each geometrically finite end we choose a smooth convex
surface as given by Theorem \ref{convex} and we let $S_0$ be the
union of these surfaces. We also choose a small horoball neighborhood
for each rank two cusp such that the boundary of each is a collection
of pairwise disjoint embedded Euclidean tori. We denote the union of
these neighborhoods and their boundary, $\cH_0$ and $H_0$,
respectively. Finally, we fix a small $\epsilon$ such that the
$\epsilon$-neighborhood of the singular locus, $\cC_\epsilon$, is a
collection of disjoint solid tori in the interior of  $N$, with boundary
$T_\epsilon$. Note that if any of $S_0$, $H_0$ or $T_\epsilon$
intersect we can choose smaller neighborhoods of each end such that
all three surfaces are disjoint.

Now let $M_0$ be the compact core of $M$ bounded by $S_0$, $H_0$
and $T_\epsilon$. Let $S_t$ and $H_t$ be distance $t$
surfaces from $S_0$ and $H_0$, respectively, and $T_t$ the
boundary of the $t$-neighborhood of the singular locus. Then we define
$M_t$ to be the compact core of $M$ bounded by $S_t$, $H_t$ and
$T_{\epsilon/(1+t)}$. 

The geometrically finite ends each define a projective structure. We
label the union of these projective structures, $\Sigma$, and the bundle of
germs of Killing fields over $\Sigma$, $E_\infty$. The surfaces $S_t$
define a map $\Pi$ from the geometrically finite ends to $\Sigma$.
By Lemma \ref{preserve automorphic}
$\Pi_*:H^1(M;E) \longrightarrow H^1(\Sigma;E_\infty)$ is an isomorphism on homology. A cohomology class $[\omega_\infty] \in H^1(\Sigma; E_\infty)$ is {\em conformal} if there exists an automorphic, conformal vector field $v$ on $\tilde{\Sigma}$ with canonical lift $s$ such that $ds \in [\omega_\infty]$.
A cohomology class $[\omega] \in H^1(M;E)$ is {\em conformal at infinity} if $\Pi_*[\omega]$
is conformal.

A cohomology class that is conformal at infinity has a representative $E$-valued 1-form that has a certain {\em standard form}. Namely there is a conformal automorphic vector field on $\tilde{\Sigma}$ with canonical lift $s_\infty$ such that $\omega = d \Pi^* s_\infty$ on the geometrically finite ends. On the rank two cusps, $\cH_0$, we can assume that $\omega$ is of the form given in \S\ref{cuspdef}. For the tubular neighborhoods of the cone singularity, $\cC_\epsilon$, standard models for $\omega$ are given on p. 36 of \cite{Hodgson:Kerckhoff:cone}. When $\omega$ is in standard form Theorem \ref{all l2} implies that $\delta \omega$, $\tr \Re \omega$ and $\hat{d} (\tr \Re \omega)$ have finite $L^2$-norm on the geometrically finite ends. The standard models on $\cH_0$ and $\cC_\epsilon$ are $d$ of the canonical lift of a divergence free, harmonic vector field. In particular, $\delta \omega$ and $\tr \omega$ are zero on $\cH_0$ and $\cC_\epsilon$. Together this implies that $\delta \omega$, $\tr \Re \omega$ and $\hat{d}( \tr \Re \omega)$ have finite $L^2$-norm on all of $M$.

We would like to show that every cohomology class in $H^1(M; E)$ is represented by a Hodge form. For example if $\tilde{\omega}$ is a closed $E$-valued and we could find a section $\tau$ of $E$ such that
\begin{equation}
\label{solvefortau}
\Delta \tau = \delta \tilde{\omega}
\end{equation}
then $\omega = \tilde{\omega} - d\tau$ we be closed and co-closed and in the same cohomology class as $\tilde{\omega}$.

To solve equation (\ref{solvefortau}) we view $\Delta$ as a linear operator on the Hilbert space of $L^2$-sections of $E$. We can solve the equation if we can show that $\Delta$ is a self-adjoint operator with positive spectrum. Since $\Delta$ is an unbounded operator we need to restrict the domain of $\Delta$.

To get a representative that is Hodge form it turns it out that we need to actually solve an equivalent equation in terms of real-valued 1-forms. 
Following \cite{Hodgson:Kerckhoff:cone} we define
$$\dom \hat{\Delta} = \{ \alpha \in L^2 | \hat{d} \alpha, \hat{\delta}
\alpha, \hat{d}\hat{\delta} \alpha, \hat{\delta} \hat{d} \alpha \in
L^2\}$$ where all derivatives are defined as distributions.
We then have the following theorem:
\begin{theorem}
\label{selfadjoint}
On a hyperbolic cone-manifold $\hat{\Delta}$ is an elliptic,
non-negative, self-adjoint operator.
\end{theorem}

{\bf Proof.} As explained in the appendix of
\cite{Hodgson:Kerckhoff:cone} the result follows from the following
Stokes' theorem. 
\begin{theorem}
\label{stokes}
Let $M$ be a hyperbolic cone-manifold. If smooth real valued forms
$\alpha$ and $\beta$ on $M$ are in $L^2$ then
\begin{equation}
\label{stokeseq}
(\hat{d}\alpha, \beta) = (\alpha, \hat{\delta}\beta).
\end{equation}
\end{theorem}

{\bf Proof.}
If $N$ is closed this is proved in
\cite{Hodgson:Kerckhoff:cone}. If $\cC$ is empty then the result is
due to Gaffney \cite{Gaffney:stokes}. More precisely Hodgson and
Kerckhoff's work shows that if $\alpha$ and $\beta$ have support on a compact
neighborhood of the singular locus then the (\ref{stokeseq})
holds. Gaffney's work shows that if the support is the
complement of a neighborhood of the singular locus than (\ref{stokeseq})
holds. General $\alpha$ and $\beta$ are the sums of forms of each type
which implies Theorem \ref{stokes} and hence Theorem
\ref{selfadjoint}. \qed{stokes}

\qed{selfadjoint}

We can now prove our Hodge theorem:
\begin{theorem}
\label{hodge}
Let $M$ be a hyperbolic cone-manifold and $\tilde{\omega}$ a smooth $E$-valued 1-form in standard form representing a cohomology class in $H^1(M;E)$ that is conformal at infinity. Then there exists a unique Hodge form $\omega$ such that following holds:
\begin{enumerate}
\item $\omega$ is cohomologous to $\tilde{\omega}$;

\item there exists an $L^2$ section $s$ of $E$ such that $ds = \tilde{\omega} - \omega$;

\item $\tilde{\omega}- \omega$ has finite $L^2$-norm on $M\backslash \cC_\epsilon$.
\end{enumerate}
\end{theorem}

{\bf Proof.} The proof is essentially the same as the proof of Theorem
2.7 in \cite{Hodgson:Kerckhoff:cone}. We begin with a brief review of
their proof and then fill in those details that their result does not
provide. 

Let $\tilde{s}$ be an automorphic section such that $\tilde{\omega} = d\tilde{s}$. Any automorphic section can be written as the sum of a canonical lift and $\imath$ times an equivariant section. In particular there exists an automorphic vector field $\tilde{v}$ and an equivariant vector field $w$ such that $\tilde{s} = \tilde{V} - \imath \curl \tilde{V} + \imath W$. Furthermore since $\tilde{\omega}$ is in standard form $w \equiv 0 $ on $\cC_{\epsilon} \cup \cH_{1}$.

By Theorem \ref{selfadjoint} there is a unique vector
field $v_1$
 solving the equation
$$ (\hat{\Delta} + 4)\hat{v}_1 = \widehat{\Delta \tilde{v}}$$
with $\hat{v}_1$ 
in $\dom \hat{\Delta}$. Since
$\hat{\Delta}$ is elliptic, $v_1$ 
is smooth. If we let $v = \tilde{v} - v_{1}$ then by \eqref{delta=deltahat+4} $v$ is harmonic so $\omega = d(V - \imath \curl V)$ is co-closed. To finish the proof we need to show that $\omega$ is a traceless and hence a Hodge form and that $\omega$ satisfies (2) and (3).

By Theorem \ref{structure} the $E$-valued 1-from $\omega$ is traceless if $\div v = 0$. Hodgson and
Kerckhoff show that
\begin{equation} 
\label{diveq}
(\hat{\Delta} + 4) \div v = 0.
\end{equation}
If we can show that $\div v$ is in $\dom \hat{\Delta}$ then we must have $\div v= 0$ since by Theorem \ref{selfadjoint} $\hat{\Delta}$ has non-negative spectrum. By construction $\div v = \div \tilde{v} - \div v_{1}$. Since $\tilde{\omega}$ is in standard form both $\div \tilde{v}$ and $\hat{d}(\div \tilde{v})$ are in $L^2$. We also know that $\hat{v}_1 \in \dom \hat{\Delta}$ so $\div v_1 = *\hat{d} \hat{v}_1$ and $\hat{d}(\div v_1)$ are in $L^2$. Together this implies that $\div v$ and $\hat{d}(\div v)$ are in $L^2$. By (\ref{diveq}) $\hat{\Delta} \div v= -4\div v$. Since $\div v$ is in $L^2$ this implies that $\hat{\Delta}\div v$ is also in $L^2$. Therefore $\div v$ is in $\dom \hat{\Delta}$ and must be zero.

We now prove (2). Let $s=  V_{1} - \imath \curl V_{1} - \imath W$. Then $ds = \tilde{\omega} - \omega$. We need to show that $s$ is in $L^{2}$ on all of $M$. First we note that $\hat{v}_{1} \in \dom \hat{\Delta}$ so $v_{1}$ and $\widehat{\curl v_{1}} = -\frac{1}{2}*\hat{d}\hat{v}_{1}$ are in $L^{2}$ on $M$. By Theorem \ref{structure}, $w= \skew \Re \tilde{\omega}$. Since $\tilde{\omega}$ is in standard form $\tilde{\omega}$ and therefore $w$ are in $L^{2}$ on the geometrically finite ends. On $\cC_{\epsilon} \cup \cH_{1}$ $w \equiv 0 $ so $w$ is in $L^{2}$ on all of $M$ proving (2).

To prove (3) we need the following lemma.

\begin{lemma}
\label{dsL2}
If $s$ is a section of $E$ such that $s$ and $\Delta s$ are in
$L^2$ on $M \backslash \cC_\epsilon$, then $ds$ is in $L^2$ on $M
\backslash \cC_\epsilon$.
\end{lemma}

{\bf Proof.} Let $f,g:M \longrightarrow [0,1]$ be smooth
functions on $M$ such that $f^2 + g^2 = 1$ and with $g =
1$ on $\cC_{\epsilon/2}$ and $g=0$ on $M \backslash \cC_\epsilon$.
Using standard techniques we can find smooth
functions $f_n :M \longrightarrow [0,1]$ such that each $f_n$ has compact
support, $|df_n|$ is bounded and $f_n \rightarrow f$ uniformally on compacts sets as $n
\rightarrow \infty$. Recall that
$$(\alpha, \beta) = \int_M \alpha \wedge *\beta^\sharp.$$
We then have
\begin{eqnarray}
(\Delta s, f^{2}_n s) & = & (d s, d(f^{2}_n s)) \nonumber \\
& = & (ds, 2f_n sdf_n) + (ds, f^{2}_n ds) \nonumber \\
& = & (f_n ds, 2sdf_n) + (f_n ds, f_n ds) \nonumber
\end{eqnarray}
where the first equality holds because $f^{2}_n s$ has compact
support. The inequality
$$\frac12(f_n ds, f_n ds) + 2(sdf_n, sdf_n) \geq |
(f_n ds, 2sdf_n) | $$
gives us
\begin{equation}
\label{ds bound}
\frac12 (f_n ds, f_n ds) \leq |(\Delta s, f^{2}_n s)| + 2(df_n s,
df_n s). 
\end{equation}
As $n \rightarrow \infty$, $(f_n ds, f_n ds) \rightarrow (fds, fds) \geq \|
ds \|_{M\backslash \cC_\epsilon}^2$ while the right hand side of
\eqref{ds bound} remains bounded since both $|f_n|$ and $|df_n|$ are
bounded for all $n$. The lemma follows.
\qed{dsL2}

To finish the proof of (3) we note that by (2) $s$ is in $L^{2}$. Since $\tilde{\omega}$ is in standard form $\Delta s = \delta \tilde{\omega}$ is also in $L^{2}$. Therefore Lemma \ref{dsL2} implies that $ds$ is in $L^{2}$ on $M \backslash \cC_{\epsilon}$ proving (3).
\qed{hodge}

{\bf Remark.} Lemma \ref{dsL2} is essentially due to Gaffney,
\cite{Gaffney:stokes}. The main difficulty is constructing the
functions, $f_n$, through a distance function which may not be
smooth. To make the functions smooth, Gaffney applies a smoothing operator to the distance
function. The convex surfaces in the geometrically finite ends allow
us to construct a smooth distance function directly.

\bigskip

A non-trivial simple closed curve $\gamma$ on $T_\epsilon$ is a
meridian if $\gamma$ is homotopically trivial in $N$. An
$E$-valued 1-form $\omega \in H^1(M;E)$ {\em preserves the cone angle}
if the infinitesimal change in holonomy of $\gamma$ induced  by $\omega$ is trivial. The asymptotic behavior of $\omega$ is described in the following
result of Hodgson and Kerckhoff \cite{Hodgson:Kerckhoff:cone}.

\begin{lemma}
\label{coneboundary} Let $M$ be a hyperbolic cone-manifold with all
cone angles $\leq 2\pi$. If $\omega \in H^1(M;E)$ is an $E$-valued 1-form
that preserves the cone angles, there exists $\epsilon_n
\rightarrow 0$ such that 
$$ \int_{T_{\epsilon_n}} \imath \omega \wedge \omega^\sharp \rightarrow
0.$$
\end{lemma}

We now prove the main theorem of this section.

\begin{theorem}
\label{vanish}
Let $M$ be a hyperbolic cone-manifold that is geometrically finite
without rank one cusps and assume that all cone angles are $\leq
2\pi$. If $\omega \in H^1(M,E)$ is an $E$-valued 1-form that is
conformal at infinity and preserves all cone angles and cusps then
$\omega \sim 0$. 
\end{theorem}

{\bf Proof.} By Theorem \ref{hodge} we can assume that $\omega$ is a Hodge form and that $\omega$ is in $L^2$ on $M \backslash \cC_\epsilon$. We will show that $\omega = 0$.

By Proposition \ref{boundaryterm}
\begin{equation}
\label{mt formula}
2\int_{M_t} \|\omega\|^2 = \int_{\del M_t} \imath \omega \wedge \omega^\sharp =B(t).
\end{equation}
We will show that $B(t) \rightarrow 0$ as $t \rightarrow
\infty$.

By Lemma \ref{coneboundary} there exists $t_i \rightarrow \infty$ such
that
$$\int_{T_{\epsilon/(1+t_i)}} \imath \omega \wedge \omega^\sharp
\rightarrow 0$$
so we are left to analyze the boundary term on $\Sigma_t \cup H_t$.

Since $\omega$ is conformal at infinity and cusp preserving,
Theorem \ref{hodge} implies that $\omega$ is in $L^2$ on $M \backslash
\cC_\epsilon$ so
\begin{eqnarray}
\int_{M\backslash \cC_\epsilon} \omega \wedge * \omega^\sharp
& = & \int_{M_0} \omega \wedge * \omega^\sharp \nonumber \\
& + & \int^\infty_0 \int_{\Sigma_t \cup H_t}
*(\omega \wedge * \omega^\sharp) dA_t dt \nonumber
\end{eqnarray}
is finite. (Note that $*(\omega \wedge * \omega^\sharp)$ is a smooth real
function since $\omega \wedge * \omega^\sharp$ is a smooth real 3-form.)
Therefore
$$ \lim_{t \rightarrow \infty} \int_{\Sigma_t \cup H_t} *(\omega \wedge *
\omega^\sharp) dA_t = 0. $$ 
We also have $|*_t( \imath \omega \wedge \omega^\sharp) | < *(\omega \wedge
* \omega^\sharp)$ where $*_t$ is the Hodge $*$-operator of the induced
metric on $\Sigma_t \cup H_t$.
Therefore
$$ \left|\int_{\Sigma_t \cup H_t} \imath \omega \wedge \omega^\sharp
\right| \leq 
\int_{\Sigma_t \cup H_t} |*_t(\imath \omega \wedge \omega^\sharp)| 
dA_t \leq \int_{\Sigma_t \cup H_t} *(\omega \wedge * \omega^\sharp) dA_t.$$
from which it follows $B(t) \rightarrow 0$ as $t \rightarrow
\infty$. Taking the limit of \eqref{mt formula} we see $2
\int_{M}\| \omega \|^2 =0$ and therefore $\omega = 0$. \qed{vanish}


\section{Representation varieties of cone-manifolds}
\label{representations}

To understand local deformations of hyperbolic structures on a
geometrically finite cone manifold we will study the representation
variety of the fundamental groups of both the manifold and its
boundary surfaces. 

Let $\Gamma$ be a finitely presented group and $G$ a Lie group. Then
$\cR(\Gamma, G)$ is the space of representations of $\Gamma$ in $G$. If
$\Gamma$ has $n$ generators and $m$ relations, $r_i$, then we can
identify $\cR(\Gamma, G)$ with a subset of $G^n$ by
$$ \cR(\Gamma, G) = \{ \gamma \in G^n : r_i(\gamma) = \id, i=1,2,\dots ,
m \}.$$
If $G$ is an algebraic group then $\cR(\Gamma, G)$ is an algebraic
variety.

We will be interested in the case where $\Gamma$ is
the fundamental group of a geometrically finite cone-manifold or a
surface with a projective structure and $G=PSL_2\cx$, the
group of hyperbolic isometries and projective transformations. For
simplicity of notation let $\cR(M) = \cR(\pi_1(M), PSL_2\cx)$ and $\cR(S) =
\cR(\pi_1(S), PSL_2\cx)$ where $M$ is a 3-manifold and $S$ a closed
surface.

The following theorem of Thurston, mentioned in the
introduction, is key to the existence of 3-dimensional hyperbolic
cone-manifolds.
\begin{theorem}[\cite{Thurston:book:GTTM}, \cite{Culler:Shalen:varieties}]
\label{lowerbound}
Let $M$ be a compact hyperbolic manifold with boundary and holonomy
representation $\rho$. Assume  that the components of $\del M$ contain
no spheres, $t$ tori and surfaces of higher genus. If $T \subset \del
M$ is a torus, we also 
assume that $\rho(\pi_1(T)) \neq 1$. Then the dimension of the
component of $\cR(M)$ containing $\rho$ is at least $t - 3\chi(M) +
3$.
\end{theorem}

Hodgson and Kerckhoff proved the following result when $\del M$ contains only tori:
\begin{theorem}[\cite{Hodgson:Kerckhoff:cone}]
\label{smooth}
Let $M$ be a compact, connected 3-manifold with non-empty boundary
consisting of $t$ tori and surfaces of higher genus.  Let $\rho \in
\cR(M)$ be an irreducible 
representation such that if $T$ is a torus component of $\del M$ then
$\rho(T) \neq 1$ or $\integers_2 \oplus \integers_2$. If the natural
map
\begin{equation}
\label{zeromap}
H^1(M, \del M; E) \longrightarrow H^1(M;E)
\end{equation}
is zero, then at $\rho$, $\cR(M)$ is a smooth complex manifold of dimension
$t - 3\chi(M) + 3$.
\end{theorem}

{\bf Sketch of proof.} To show that a variety is smooth one needs to
show that the dimension of the Zariski tangent space is
minimal. Theorem \ref{lowerbound} gives a lower bound for this
dimension so we need to show that the dimension at $\rho$ equals this
lower bound.

A fundamental result of Weil shows that $\dim T\cR(M)_\rho = H^1(M;E) +
3$ if $\rho$ is irreducible. Hodgson and Kerckhoff show that if the
natural map (\ref{zeromap}) is zero then 
$$\dim H^1(M;E) = \frac12 \dim H^1(\del M;E(\del M)).$$
We are left to calculate $\dim H^1(\del M;E)$ which will be the sum of
the dimensions of $H^1(S;E)$ at $\rho$ over all connected components
$S$ of $\del M$. The dimension of $H^1(S;E)$ is well known. For a
torus $T$ with representation $\rho(T) \neq 1$ or $\integers_2 \oplus
\integers_2$, $\dim H^1(T;E) = 2$. If $S$ has genus $>1$ at an
irreducible representation, $\dim H^1(S; E) = -3\chi(S)$. Summing these
dimensions we find
$$\dim H^1(M;E) = \frac12 \dim H^1(\del M;E) = \frac12(2t - 3\chi(\del
N)) = t - 3\chi(M).$$
Since the dimension of the tangent space at $\rho$ is minimal,
$\cR(M)$ is smooth and has dimension $t - 3\chi(M)$. \qed{smooth}

{\bf Remark.} To turn our sketch into an actual proof we need to view
$\cR(M)$ as a scheme instead of a variety. Then the Zariski tangent
space of $\cR(M)$ is the space of 1-cocycles with coefficients
in the module $\Ad \rho$. Furthermore at a representation $\rho$ satisfying the conditions of Theorem \ref{smooth} the algebraically defined Mumford quotient $\cR(M)//PSL_2\cx$ is isomorphic to the topological quotient $\cR(M)/PSL_2\cx = R(M)$. At the image of $\rho$, $R(M)$ will be a complex manifold whose (differentiable) tangent space is canonically identified with $H^1(M; E)$.

\bigskip
%
%
%
%
To apply this result to geometrically finite hyperbolic cone manifolds
we need the following result:
\begin{prop}
\label{irreducible}
Let $M$ be a hyperbolic cone-manifold that is geometrically finite
without rank one cusps and let $\rho$ be its holonomy representation.
\begin{enumerate}
\item \label{irredend} The restriction of $\rho$ to each geometrically
finite end is irreducible.

\item \label{rho irreducible} $\rho$ is irreducible.

\item \label{rhot non-trivial} Let $T$ be the boundary of an
$\epsilon$-neighborhood of a 
component of the singular locus. Then the image of $\rho(\pi_1(T))$ is
infinite and non-parabolic.
\end{enumerate}
\end{prop}

{\bf Proof.} \ref{irredend}. The holonomy of a projective
structure on a surface of genus $> 1$ is always irreducible, for a
reducible representation fixes a point on $\chat$ and hence has image
an affine group. Since a surface of genus $>1$ cannot have an affine
structure this is impossible. The restriction of $\rho$ to a
geometrically finite end is also the holonomy of a projective
structure and therefore is irreducible.

\ref{rho irreducible}. If $\vol(M)$ is finite
then this is Lemma 4.6 in \cite{Hodgson:Kerckhoff:cone}. If not $M$
contains a geometrically finite end on which by (\ref{irredend}) the
holonomy is irreducible which implies that $\rho$ is irreducible.

\ref{rhot non-trivial}. The holonomy of any homotopically
non-trivial closed curve on $T_\epsilon$ that is not a multiple of the
meridian will have hyperbolic holonomy. This implies that the image of
$\rho(\pi_1(T))$ is infinite and non-parabolic. \qed{irreducible}

The following corollary follows directly from Theorems \ref{vanish}
and \ref{smooth} along with Proposition \ref{irreducible}.
\begin{cor}
\label{smoothcone}
Let $\rho$ be the holonomy representation of a hyperbolic
cone-manifold $M$ that is geometrically finite without rank one
cusps. If all cone angles of $M$ are $\leq 2 \pi$ then $R(M)$ is
smooth at $\rho$ with dimension $n+m - 3\chi(M)$ where $n$ is the
number of components of $\cC$ and $m$ is the number of rank two cusps.
\end{cor}



For $\gamma \in \pi_1(M)$ let $\cL_\gamma (\rho)$ denote the complex length of $\rho(\gamma)$. If $\rho(\gamma)$ is hyperbolic then $\cL_\gamma(\rho)$ is the sum of the translation length plus $\imath$ times the angle of rotation. While this is only well defined up to sign and the angle is only defined modulo $2\pi$ after making an initial choice $\cL_\gamma$ extends to a holomorphic function in a neighborhood of $\rho$. In our setting when $\gamma$ is the meridian of a cone singularity it is natural to choose $\cL_\gamma(\rho)$ to be the cone angle. If $\rho(\gamma)$ is parabolic we define $\cL_\gamma(\rho)=0$. In this case there is no way to make a choice of sign. Instead we view $\cL_\gamma$ as a map to $\cx /\{\pm 1\}$. Although this will allow us to extend $\cL_\gamma$ to a continuous map in neighborhood of $\rho$ it will not in general be differentiable. For this reason at parabolic elements it is convenient to use the trace map. That is $\Tr_\gamma(\rho)$ assigns to each $\rho \in R(M)$ the trace of $\rho(\gamma)$. Again this map is only defined up to sign but at a parabolic the trace is $\pm 2$ so a well defined choice of sign can be made. The trace then extends to a holomorphic map at parabolics. Note that $\Tr_\gamma(\rho) = 2 \cosh(\cL_\gamma(\rho)/2)$.

To understand the derivative of $\cL_\gamma$ (or $\Tr_\gamma$) it is helpful to look at the bundle $E(\gamma)$ which we define to be the restriction of $E$ to a smooth loop in the free homotopy class of $\gamma$. Then each cohomology class $\omega \in H^1(M;E)$ restricts to a cohomology class in $H^1(\gamma; E(\gamma))$. If $\rho(\gamma)$ is not the identity then $H^1(\gamma; E(\gamma)) \isom \cx$ where the natural isomorphism sends cohomology classes to tangent vectors to the space of complex lengths (or traces). Note that the infinitesimal change in holonomy of $\gamma$ induced by $\omega$ is trivial if and only if $\omega$ restricts to a trivial element of $H^1(\gamma; E(\gamma))$.
More precisely we have the following lemma which is essentially contained in Theorem 4.5 in \cite{Hodgson:Kerckhoff:cone}.
%
%

\begin{lemma}
\label{Lw=0}
\begin{enumerate}
\item Let $\gamma$ be a meridian of the cone singularity. Then $(\cL_\gamma)_* \omega = 0$ if and only if $\omega$ preserves the cone angle.

\item Let $\gamma$ be homotopic to a rank two cusp. Then $(\Tr_\gamma)_* \omega = 0$ if and only if $\omega$ is cusp preserving.
\end{enumerate}
\end{lemma}

Note that if the cone angle $2\pi$ the holonomy of the meridian will be the identity. This special case is also dealt with in Theorem 4.5 of \cite{Hodgson:Kerckhoff:cone}.

%

We now describe a local parameterization of $R(M)$ that is the main
theorem of this paper. To do so we need to recall some basic facts
about the space of marked projective structures, $P(S)$, on a closed
surface $S$ of genus $>1$. These can all be found in
{\cite{Gunning:book:LVB}. $P(S)$ is a complex manifold of dimension
$-3\chi(S)$. If $\Sigma \in P(S)$ is a projective structure then the
tangent space of $P(S)$ at $\Sigma$ can be canonically identified with
$H^1(\Sigma;E(\Sigma))$. The {\em Teichm\"uller space}, $T(S)$, of $\Sigma$ is the
space of marked conformal structures on $S$. Since a projective
structure also defines a conformal structure there is a projection, $p:P(S)
\longrightarrow T(S)$. Furthermore, if $\omega \in H^1(\Sigma; E(\Sigma))$ is an
$E(\Sigma)$-valued 1-form than $p_*\omega = 0$ if and only if $\omega$ is
conformal. There is also a holonomy map, $h:P(S) \longrightarrow
R(S)$. We will need the following theorem:
\begin{theorem}[Hejhal \cite{Hejhal:monodromy}]
\label{Hejhal}
The map $h$ is a holomorphic, local homeomorphism.
\end{theorem}

Assume $M$ is a hyperbolic cone-manifold that is geometrically finite
without rank one cusps. Assume the cone singularity has $n$
components and that $M$ has $m$ rank two cusps. Let $S$ be the union
of the higher genus boundary components of $\del M$. $R(S)$, $P(S)$ and
$T(S)$ will be the product of the representation varieties, spaces of
projective structures and Teichm\"uller spaces, respectively, of the
components of $S$. For each component of
the cone singularity we let $\cL_i$, $i=1, \dots n$, be the complex length of the meridian. For each rank two cusp we choose a generator of the corresponding $\integers \oplus \integers$ subgroup and let $\cL_i$ be its complex length and $\Tr_i$ its trace with $i= n+1, \dots n+m$. We also have a map,
$\del:R(M) \longrightarrow R(S)$, that restricts each representation
to a representation of the boundary surfaces. We then define a maps,
$\Phi$ and $\overline{\Phi}$
by
$$\Phi(\sigma) = (\cL_1(\sigma), \dots, \cL_{n}(\sigma), \Tr_{n+1}(\sigma), \dots, \Tr_{n+m}(\sigma), p \circ
h^{-1} \circ \del (\sigma))$$
and
$$\overline{\Phi}(\sigma) = (\cL_1(\sigma), \dots, \cL_{n+m}(\sigma), p \circ h^{-1} \circ \del(\sigma))$$
for $\sigma \in R(M)$. Note that the image of $\Phi$ is contained in $\cx^{n+m} \times T(S)$ while the image of $\overline{\Phi}$ is contained in $\cx^n \times (\cx/\pm 1)^m \times T(S)$.

We now prove our main theorem.
\begin{theorem}
\label{parameterization}
Assume $M$ is a geometrically finite hyperbolic cone-manifold without rank one cusps and with holonomy representation
$\rho$. If all cone angles of $M$ are $\leq 2\pi$ then
$\Phi$ is a holomorphic, local homeomorphism and $\overline{\Phi}$ is a local homeomorphism at $\rho$.
\end{theorem}

{\bf Proof.} By Corollary \ref{smoothcone}, $R(M)$ is a smooth complex
manifold of dimension $n+m - 3\chi(M)$ which is equal to the dimension
of $\cx^{n+m} \times T(S)$. Since $\cL_i$, $\Tr_i$ and $p \circ h^{-1}$ are all holomorphic $\Phi$ is holomorphic. To show that $\Phi$ is
a local homeomorphism we need to show that $\Phi_*$ has trivial
kernel. If $\omega \in H^1(M;E)$ is an $E$-valued 1-form such that
$\Phi_* \omega = 0$ then $(\cL_i)_*\omega=0$, $(\Tr_i)_* \omega = 0$ and $p_* \omega=0$. By Lemma \ref{Lw=0},
$(\cL_i)_* \omega = 0$ implies that $\omega$ preserves the
cone angle and $*\Tr_\gamma)_*\omega = 0$ implies that
$\omega$ is cusp preserving.  Finally if $p_* \omega =0$, then
$\omega$ is conformal at infinity. Therefore Theorem \ref{vanish}
implies that $\omega$ is trivial so $\Phi_*$ has trivial kernel and
$\Phi$ is a local homeomorphism at $\rho$. The relationship between the trace and the complex length then implies that $\overline{\Phi}$ is a local homeomorphism. \qed{parameterization}

This parameterization leads to our local rigidity theorem.

\begin{theorem}
\label{main}
If $M$ is a geometrically finite cone-manifold without rank one cusps
and all cone angles are $\leq 2 \pi$ then $M$ is locally rigid rel
cone angles and the conformal boundary.
\end{theorem}

{\bf Proof.} Let $M_t$ be a smooth family of cone-metrics on $M$ such
that $M_0 = M$ and such that the conformal structures at infinity and
cone angles of $M_t$ agree with those of $M$. Then by Theorem
\ref{parameterization} the holonomy representations $\rho_t$ for $M_t$
are equal to $\rho_0$. 

Theorem 1.7.1 of \cite{Canary:Epstein:Green} implies that for every
compact core, $M'$, of
$M_0$ there exists a $t'$ such that $M'$ isometrically embeds in $M_t$
for $t < t'$. Choose $M'$ such that $\del M'$ is a collection of
convex surfaces of higher genus and Euclidean tori around each rank
two cusp and component of the singular locus. Then any isometry of
$M'$ into $M_t$ can be extended to an isometry from $M_0$ onto
$M_t$. Hence $M_0$ is locally rigid rel cone angles and the conformal
boundary. \qed{main}

\bibliographystyle{../tex/math}
\bibliography{../tex/math}
\end{document}

%% file: macros.tex
\newcommand{\bb}{\mathbb}

\newcommand{\cx}{{\bb C}}

\newcommand{\integers}{{\bb Z}}

\newcommand{\reals}{{\bb R}}

\newcommand{\hthree}{{\bb H}^3}

\newcommand{\ddz}{\frac{\partial}{\partial z}}
\newcommand{\qed}[1]{\nopagebreak[4]{\tiny \hfill
\fbox{\ref{#1}} \linebreak }\pagebreak[2]}
\newcommand{\del}{\partial}

\newcommand{\isom}{\cong}

\newcommand{\zbar}{{\overline{z}}}

\newcommand{\chat}{\widehat{\cx}}

\newcommand{\ad}{\operatorname{ad}}
\newcommand{\Ad}{\operatorname{Ad}}

\newcommand{\area}{\operatorname{area}}

\newcommand{\ave}{\operatorname{ave}}

\renewcommand{\div}{\operatorname{div}}
\newcommand{\dom}{\operatorname{dom}}

\newcommand{\Hom}{\operatorname{Hom}}

\newcommand{\id}{\operatorname{id}}

\renewcommand{\Im}{\operatorname{Im}}
\newcommand{\Isom}{\operatorname{Isom}}

\renewcommand{\Re}{\operatorname{Re}}

\renewcommand{\skew}{\operatorname{skew}}

\newcommand{\str}{\operatorname{str}}

\newcommand{\sym}{\operatorname{sym}}

\newcommand{\Tr}{\operatorname{Tr}}
\newcommand{\tr}{\operatorname{tr}}
\newcommand{\vol}{\operatorname{vol}}
\newcommand{\curl}{\operatorname{curl}}
\newcommand{\ex}{\operatorname{ex}}

\newtheorem{theorem}{Theorem}[section]
\newtheorem{prop}[theorem]{Proposition}
\newtheorem{lemma}[theorem]{Lemma}
\newtheorem{cor}[theorem]{Corollary}

\newtheorem{defn}[theorem]{Definition}

\newcommand{\cC}{{\cal C}}

\newcommand{\cH}{{\cal H}}

\newcommand{\cL}{{\cal L}}

\newcommand{\cR}{{\cal R}}
